\documentclass[reqno,11pt]{amsart}
\usepackage{etex}
\usepackage{tikz}
\usetikzlibrary{matrix,arrows}

\usepackage{amsmath,amssymb,amscd,amsxtra,dsfont}
\usepackage{graphicx,pst-grad,pst-plot,pst-coil}
\usepackage{pstricks}
\usepackage{eucal}
\usepackage[margin=1.5in]{geometry}

\usepackage[colorlinks=true,linkcolor=blue,citecolor=blue]{hyperref}  

\usepackage{mathptmx}
\DeclareSymbolFont{txsymbols}{OMS}{txsy}{m}{n}
\SetSymbolFont{txsymbols}{bold}{OMS}{txsy}{bx}{n}

\DeclareSymbolFont{txlargesymbols}{OMX}{txex}{m}{n}
\SetSymbolFont{txlargesymbols}{bold}{OMX}{txex}{bx}{n}

\let\amalg\relax
\DeclareMathSymbol{\amalg}{\mathbin}{txsymbols}{113}
\let\coprod\relax
\DeclareMathSymbol{\coprod}{\mathop}{txlargesymbols}{96}
\DeclareMathSymbol{\nabla}{\mathord}{txsymbols}{114}

\input diagxy
\input xy
\xyoption{all}
\def\bfig{\vcenter\bgroup\xy}
\def\efig{\endxy\egroup}

\theoremstyle{plain}

\newtheorem{thm}{Theorem}
\newtheorem{proposition}[thm]{Proposition}
\newtheorem{corollary}[thm]{Corollary}
\newtheorem{lemma}[thm]{Lemma}

\theoremstyle{definition}
\newtheorem{remark}[thm]{Remark}

\newenvironment{definition}{\refstepcounter{thm}
{\medskip\par\noindent\bf Definition \arabic{section}.\arabic{thm}.
}}{\vskip 2ex\par}
\newenvironment{example}{\refstepcounter{thm}
{\medskip\par\noindent\it Example \arabic{section}.\arabic{thm}.}}{\vskip 2ex\par}
\newenvironment{exam}{\refstepcounter{thm}
{\medskip\par\noindent\it Example \arabic{section}.\arabic{thm}.}}{\vskip 2ex\par}
\newenvironment{example*}{
{\medskip\noindent\it Example.}}{\vskip 2ex\par}

\newenvironment{exer}{\refstepcounter{thm}
{\medskip\par\noindent\sc Exercise \arabic{section}.\arabic{thm}
}}{\vskip 2ex\par}

\renewenvironment{proof}{
{\medskip\par\noindent\sc Proof.\ }}{\vskip 2ex\par}

\newcommand{\ga}{\alpha}

\newcommand{\cala}{\mathcal{A}}

\newcommand{\calc}{\mathcal{C}}

\newcommand{\cale}{\mathcal{E}}
\newcommand{\calf}{\mathcal{F}}
\newcommand{\calg}{\mathcal{G}}
\newcommand{\calh}{\mathcal{H}}

\newcommand{\call}{\mathcal{L}}

\newcommand{\calo}{\mathcal{O}}
\newcommand{\calq}{\mathcal{Q}}

\newcommand{\cals}{\mathcal{S}}

\newcommand{\calz}{\mathcal{Z}}

\newcommand{\R}{\mathbb{R}}

\newcommand{\comp}{\mathrel{\scriptstyle\circ}}

\newcommand{\C}{\mathbb{C}}

\newcommand{\cinf}{{C^{\infty}}}

\newcommand{\Gr}{\operatorname{Gr}}

\newcommand{\hh}{\mathbb{H}}

\newcommand{\Hom}{\operatorname{Hom}}

\newcommand{\id}{{\mathds{1}}}
\newcommand{\im}{\operatorname{im}}

\newcommand{\Open}{\operatorname{Open}}

\newcommand{\supp}{\operatorname{supp}}

\newcommand{\Z}{\mathbb{Z}}

\def \pf{\begin{proof}}
\def \epf{\end{proof}}
\def \enum{\begin{enumerate}}
\def \eenum{\end{enumerate}}

\renewcommand{\to}{\rightarrow}

\def \defi{\begin{definition}}
\def \edefi{\end{definition}}

\def \prop{\begin{proposition}}
\def \eprop{\end{proposition}}

\def \lem{\begin{lemma}}
\def \elem{\end{lemma}}

\def \cor{\begin{corollary}}
\def \ecor{\end{corollary}}

\newcommand{\ex}{\begin{example*}}
\newcommand{\eex}{\end{example*}}

\newenvironment{exer*}
  {\small\begin{exercise}}
  {\end{exercise}}

\newcommand{\probref}[1]{\textbf{\ref{#1}} } 

\def \ex*{\begin{example*}}
\def \eex*{\end{example*}}

\newenvironment{remark*}{
{\medskip\noindent\it Remark.}}{\vskip 2ex\par}
\def \rem*{\begin{remark*}}
\def \erem*{\end{remark*}}

\newenvironment{claim*}{
{\medskip\noindent\it Claim.}}{\vskip 2ex\par}

\def \pf{\begin{proof}}
\def \epf{\end{proof}}

\def \enum{\begin{enumerate}}
\def \eenum{\end{enumerate}}

\numberwithin{equation}{section}
\numberwithin{figure}{section}
\numberwithin{thm}{section}
\def\H{\mathbb H}

\def\H{\mathbb H}
\def\I{{\rm i}}
\def\C{\mathbb C}

\def\R{{\mathbb R}}
\def\Z{\mathbb Z}

\newcommand{\term}[1]{\textbf{\textit{#1}}}
 
\usepackage[colorlinks=true,linkcolor=blue]{hyperref}

\begin{document}
\title
	{Introduction to Sheaf Cohomology}
\begin{center}	
Lecture Notes for the Summer School
on Algebraic Geometry and Number Theory\\
Galatasaray University, Istanbul, Turkey\\
June 2--13, 2014\\
\end{center}

\author{Loring W. Tu}
\address{Department of Mathematics\\
Tufts University\\
Medford, MA 02155} 
\email{loring.tu@tufts.edu}
\thanks{Version 6, June 4, 2022.  I am indebted to George Leger and especially to Jeffrey D.~Carlson for their detailed comments and feedback.
}

\maketitle

\vspace{-.8cm}
\tableofcontents

\setcounter{page}{1} \setcounter{thm}{0}


This is a more or less verbatim account of my lectures at
the CIMPA Summer School on Algebraic Geometry and
Number Theory in Istanbul, Turkey, in June 2014.
My goal was to introduce to the uninitiated, in just four lectures, the wonderful
techniques of sheaf cohomology, hypercohomology, and 
spectral sequences.
Because of the diversity of the background of the audience,
I decided to start at the beginning and to assume no knowledge
of sheaves.
Indeed, since there exists a multitude of approaches to sheaves,
it may be desirable to pick one out that seems to me the best
because of its simplicity.

As prerequisites, I assume that the students have a good
knowledge of manifolds, including the exterior calculus
of differential forms, as in \cite{tu_m}.
Our main reference for sheaf cohomology is Part I of the article ``From
sheaf cohomology to the algebraic de Rham theorem''  \cite{elzein--tu}
that Fouad El Zein and I wrote.
However, because I had only four lectures at the CIMPA School,
I tried to move at a fast clip, omitting some of the proofs.
Most of the omitted proofs can be found in \cite{elzein--tu}.

Introduced in 1946 by Jean Leray and further developed 
in subsequent years by Henri Cartan (\cite{mccleary} and \cite{miller}),
sheaves are a powerful tool for relating
local and global phenomena on a space.
Smooth differential forms of a given degree define a sheaf on a manifold.
On a complex manifold, in addition to smooth differential forms,
there are holomorphic differential forms that also define
sheaves on the manifold.  

Sheaf cohomology may be viewed as a generalization of singular
cohomology from constant coefficients to variable coefficients.
Sheaf cohomology has become absolutely essential in modern algebraic geometry
as well as certain areas of topology and complex analysis.
For example, cohomology with coefficients in the sheaf $\Omega^p$
of holomorphic $p$-forms on a complex manifold
gives invariants of the complex structure.

Hypercohomology is a generalization of sheaf cohomology from
one sheaf to a complex of sheaves.  It is a functor from the category
of complexes of sheaves to the category of abelian groups.
There are two spectral sequences that converge to hypercohomology.
As we will see, the isomorphism of the limits of the two spectral sequences often
yields important isomorphism theorems in complex and algebraic geometry.

On a smooth manifold $M$, a smooth $k$-form $\omega$ can be 
integrated over a continuous $k$-chain to yield a real number.
Thus, $\int_{(\ )}\omega$ is a $k$-cochain on $M$.
Let $\cala^k(M)$ be the group of $\cinf$ $k$-forms on $M$
and $S^k(M,\R)$ the group of continuous real singular $k$-cochains
on $M$.
By Stokes' theorem, the map
\begin{align*}
\varphi\colon \cala^k(M) &\to S^k(M,\R),\\
\omega &\mapsto \int_{(\ )}\omega,
\end{align*}
satisfies
\[
\varphi(d\omega) = \int_{(\ )} d\omega = \int_{\partial(\ )}\omega =
\delta \int_{(\ )} \omega = \delta (\varphi \omega),
\]
where $\delta$ is the coboundary operator on $S^k(M,\R)$.
Thus, $\varphi$ is a cochain map and induces a linear map
\[
\varphi^*\colon H_{\rm{dR}}^k(M) \to H_{\rm{sing}}^k(M,\R)
\]
from de Rham cohomology $H_{\rm{dR}}^k(M)$ to singular cohomology
$H_{\rm{sing}}(M,\R)$.
The de Rham theorem states that the induced map $\varphi^*$
in cohomology is an isomorphism.
We will show how to use sheaf cohomology, hypercohomology, and
spectral sequences to prove the de Rham theorem.
An immediate consequence of the de Rham theorem is that
although de Rham cohomology is defined in terms of smooth forms,
it is actually a topological invariant, since singular cohomology
is defined in the continuous category.

A differential form $\omega$ on a manifold $M$ is \term{smooth}
if locally it can be written as $\omega = \sum a_I dx^{i_1} \wedge
\cdots \wedge dx^{i_k}$, where $(U, x^1, \ldots, x^n)$ is a chart
and the $a_I$ are $\cinf$ functions on $U$.
On a complex manifold, one can define similarly \term{holomorphic
forms}: a differential form $\omega$ on a complex manifold
is \term{holomorphic} if locally it can be written as
$\omega = \sum a_I dz^{i_1} \wedge \cdots \wedge dz^{i_k}$,
where $(U, z^1, \ldots, z^n)$ is a complex chart and the $a_I$
are holomorphic functions on $U$.

Let $\Omega^k(M)$ be the group of holomorphic $k$-forms
on $M$.
The holomorphic forms on $M$ constitute a \term{differential
complex}, called the \term{holomorphic de Rham complex},
\[
0 \to \Omega^0(M) \overset{d}{\to} \Omega^1(M) \overset{d}{\to}
\cdots \to \Omega^n(M) \to 0,
\]
in the sense that $d \comp d = 0$.
One might ask if the singular cohomology of a complex manifold
with complex coefficients can be computed using only 
holomorphic forms.
According to the analytic de Rham theorem, the answer is yes,
but it is not as simple as taking the cohomology of the
holomorphic de Rham complex.
Here the hypercohomology of the complex $\Omega^{\bullet}$
of sheaves of holomorphic forms on $M$ makes its
appearance:
\[
H_{\rm{sing}}^*(M,\C) \simeq \mathbb{H}^*(M, \Omega^{\bullet}).
\]
In case $M$ is a Stein manifold (a complex manifold that is
a closed subset of $\C^N$), then a simpler form of analytic
de Rham theorem holds:
\[
H_{\rm{sing}}^*(M,\C) \simeq h^*\big( \Omega^{\bullet}(M)\big) :=
\text{cohomology of the holomorphic de Rham complex}.
\]

In the algebraic category, a differential $k$-form $\omega$ on a
smooth algebraic variety is \term{algebraic} if locally it can be
written as $\omega = \sum f_I dg^{i_1} \wedge \cdots \wedge
dg^{i_k}$, where the $f_I$ and $g^j$ are regular functions.
If $X$ is a smooth complex algebraic variety with the Zariski
topology, let $X(\C)$ be the set of complex points with the complex
topology.
Grothendieck's algebraic de Rham theorem states that the
singular cohomology $H_{\rm{sing}}^*\big( X(\C), \C \big)$ with
complex coefficients can be computed from the algebraic differential forms
on $X$:
\[
H_{\rm{sing}}^* \big( X(\C), \C \big) \simeq \mathbb{H}^*(X, \Omega_{\rm{alg}}^{\bullet}),
\]
where $\Omega_{\rm{alg}}^{\bullet}$ is the complex of sheaves of 
algebraic forms on $X$.
For a smooth affine variety $X$, the singular cohomology $H_{\rm{sing}}^*(M, \C)$
is the cohomology of the algebraic de Rham complex:
\[
H_{\rm{sing}}^*\big( X(\C), \C \big) \simeq 
h^*\big( \Omega_{\rm{alg}}^{\bullet}(X) \big).
\]

In this short article we will not be able to prove the algebraic de Rham theorem;
for a proof, see \cite{elzein--tu}.  We will give an introduction to sheaf cohomology,
hypercohomology, and spectral sequences, enough to prove the de Rham theorem
for smooth manifolds and the analytic de Rham theorem for complex
manifolds.

\bigskip
\bigskip
\noindent
\section*{{\bf Lecture 1. Presheaves and Sheaves}}

\section{Presheaves}

The functor $\cala^*$ that assigns to every open set $U$ on a
manifold the vector space of $\cinf$ forms on $U$ is an example of
a presheaf. By definition a
\term{presheaf}\index{presheaf} $\calf$ of abelian groups on a topological space $X$
is a function that assigns to every open set $U$ in $X$ an abelian
group $\calf(U)$ and to every inclusion of open sets $i_U^V\colon
V \rightarrow U$ a group homomorphism $\calf(i_U^V):=\rho_V^U$, called the
\term{restriction} from $U$ to $V$,
\[
\rho_V^U\colon \calf(U) \rightarrow \calf(V),
\]
such that the system of restrictions $\rho_V^U$
satisfies the following properties: 
\enum 
\item[(i)] (identity)
$\rho_U^U = \id_{F(U)}$, the identity map on $\calf(U)$; 
\item[(ii)]
(transitivity) if $W \subset V \subset U$, then $\rho_W^V \comp
\rho_V^U = \rho_W^U$. \eenum We refer to elements of $\calf (U)$
as \term{sections} of $\calf$ over $U$.
The group $\calf(U)$ is also written $\Gamma(U,\calf)$.
Elements of $\Gamma(X, \calf)$ are called \term{global sections} of $\calf$.

If $\calf$ and $\calg$ are presheaves on $X$, a
\term{morphism}\index{morphism!of presheaves} $f\colon \calf
\rightarrow \calg$ of presheaves is a collection of group
homomorphisms $f_U\colon \calf(U) \rightarrow \calg(U)$, 
indexed by open sets $U$ in $X$, that commute with the restrictions,
i.e., such that each diagram
\begin{equation} \label{6e:restrictions}
\bfig \xymatrix{ \calf(U) \ar[r]^{f_U} \ar[d]_{\rho_V^U} &
\calg(U) \ar[d]^{\rho_V^U} \\
\calf(V) \ar[r]_{f_V} &\calg(V)} \efig
\end{equation}
is commutative.
Although we write both vertical maps as $\rho_V^U$, they are in
fact not the same map; the first one is $\calf(i_U^V)$ and
the second is $\calg(i_U^V)$.
If we write $\omega|_U$ for $\rho_V^U(\omega)$, then the diagram 
\eqref{6e:restrictions} is equivalent to $f_V(\omega|_V) =
\left.f_U(\omega)\right|_V$ for all $\omega\in \calf(U)$.
In practice, we often omit the subscripts in $f_U$
and $f_V$ and write them simply as $f$.

For any topological space $X$, let $\Open(X)$ be the category in
which the objects are open subsets of $X$ and for any two open
subsets $U$, $V$ of $X$, the set of morphisms from $V$ to $U$ is
\[
\Hom(V,U) :=
\begin{cases}
\{ \text{inclusion } i_U^V\colon V \rightarrow U \} &\text{if } V \subset U,\\
\text{the empty set } \varnothing &\text{otherwise.}
\end{cases}
\]
In functorial language, a presheaf of abelian groups is simply a contravariant
functor from the category $\Open(X)$ to the category of abelian
groups, and a morphism of presheaves is a natural
transformation\index{natural transformation} from the functor
$\calf$ to the functor $\calg$. What we have defined are
presheaves of abelian groups; it is possible to define similarly
presheaves of vector spaces, of algebras, and indeed of objects in any
category, but all the presheaves that we consider will be
presheaves of abelian groups.

\begin{exam} \label{exam:zero}
The \term{zero presheaf}\index{zero presheaf} $\calf$ on a
topological space $X$ associates to every open set $U$ the zero
group $\calf(U) =0$ and to every inclusion $V \subset U$ the zero
map $\calf(U) \to \calf(V)$.
\end{exam}

\begin{example*}
If $G$ is an abelian group, we define the \term{presheaf of
locally constant $G$-valued functions on $X$}\index{presheaf!of
locally constant functions} (constant on connected components)
 to be the presheaf $\underline{G}$
that associates to every open set $U$ in $X$ the group
\[
\underline{G}(U) :=\{ \text{locally constant functions }f\colon U
\rightarrow G\}
\]
and to every inclusion of open sets $V \subset U$ the restriction
$\rho_V^U\colon \underline{G}(U) \rightarrow \underline{G}(V)$ of
locally constant functions.
\end{example*}

\section{The Stalk of a Presheaf}

On a smooth manifold $M$, the function that assigns to every 
open set $U \subset M$ the group $\cinf(U)$ of $\cinf$ real-valued functions
on $U$ is a presheaf.
As we know from manifold theory, the behavior
of $\cinf$ functions at a point is encoded in the
the \emph{germs} of the functions at the point.
The corresponding notion for a presheaf is the
\emph{stalk} of the presheaf at a point.
To define the stalk, we recall
an algebraic construction called the
\emph{direct limit} of a direct system of groups.

A \term{directed set}\index{directed set} is a set $I$ with a
binary relation $\le$ satisfying 
\enum \item[(i)] (reflexivity) for all $a \in I$,
$a \le a$,
\item[(ii)] (transitivity) for all $a, b, c \in I$, if $a \le b$ and $b
\le c$, then $a \le c$,
\item[(iii)] (upper bound) for all $a,b \in
I$, there is an element $c \in I$, called an \term{upper
bound}\index{upper bound} of $a$ and $b$, such that $a \le c$ and $b
\le c$. \eenum We often write $b \ge a$ if $a \le b$.

%
%

A \term{direct system of groups}\index{direct system of groups} is
a collection of groups $\{ G_i \}_{i \in I}$ indexed by a directed
set $I$ and a collection of group
homomorphisms $f_b^a\colon G_a \rightarrow G_b$ indexed by pairs $a \le
b$ in $I$ such that
\enum 
\item[(i)] $f_a^a = \id_{G_a}$, the identity map on $G_a$,
\item[(ii)]
$f_c^a = f_c^b \comp f_b^a$ for $a \le b \le c$ in $I$. 
\eenum On the disjoint
union $\coprod_i G_i$ we introduce an equivalence relation
$\sim$ by decreeing two elements $g_a$ in $G_a$ and $g_b$ in $G_b$
to be equivalent if there exists an upper bound $c$ of $a$ and $b$
such that $f_c^a(g_a) = f_c^b(g_b)$ in $G_c$. The \term{direct
limit}\index{direct limit} of the direct system, denoted by
$\varinjlim_{i\in I} G_i$, is the quotient of the disjoint union
$\coprod_i G_i$ by the equivalence relation $\sim$; in other
words, two elements of $\coprod_i G_i$ represent the same element
in the direct limit if they are ``eventually equal.'' We make the
direct limit $\varinjlim G_i$ into a group by defining $[g_a] +
[g_b] = [f_c^a(g_a) + f_c^b(g_b)]$, where $c$ is an upper bound of
$a$ and $b$ and $[g_a]$ is the equivalence class of $g_a$. It is
easy to check that the addition $+$ is well defined and that with this
operation the direct limit $\varinjlim G_i$ becomes a
group; moreover, if all the groups $G_i$ are abelian, then so is their
direct limit. Instead of groups, one can obviously also consider
direct systems of modules, rings, algebras, and so on.

\begin{example*}
Fix a point $p$ in a manifold $M$ and let $I$ be the directed set
consisting of all neighborhoods of $p$ in $M$, with $\le$ being
reverse inclusion: $U \le V$ if and only if $V \subset U$. Let
$\cinf(U)$ be the ring of $\cinf$ functions on $U$. Then $\{
\cinf(U) \}_{U \ni p}$ is a direct system of rings and its direct
limit $C_p^{\infty} := \varinjlim_{U \ni p} \cinf(U)$ is precisely the ring of
germs of $\cinf$ functions at $p$.
\end{example*}

If $\calf$ is a presheaf of abelian
groups on a topological space $X$ and $p$ is a point in $X$, then
$\{ \calf(U)\}_{U \ni p}$, where $U$ ranges over all
open neighborhoods of $p$, is a direct system of abelian groups.
The direct limit $\calf_p := \varinjlim_{U \ni p} \calf(U)$ is
called the \term{stalk} of $\calf$ at $p$. An element of the stalk
$\calf_p$ is called a \term{germ} of sections at $p$.
For example, the ring $C_p^{\infty}$ is the stalk at $p$ of the
presheaf $\cinf(\ )$ of $\cinf$ functions on the manifold $M$.

A morphism of presheaves $\varphi\colon \calf \to \calg$ over a
topological space $X$ induces a
morphism of stalks $\varphi_p\colon \calf_p \to \calg_p$ at each $p
\in X$ by
sending the germ at $p$ of a section $s \in \calf(U)$ to the germ at $p$ of the
section $\varphi(s) \in \calg(U)$.
The morphism $\varphi_p\colon \calf_p \to \calg_p$ of stalks
is also called the \term{stalk map} at $p$.

\section{Sheaves}

The stalk of a presheaf at a point embodies in it the local
character of the presheaf about the point. However, in general
there is no relation between the global sections and the stalks of
a presheaf. 

\begin{exam}\label{exam:nonzero}
If $G$ is an abelian group and $\cals$ is
the presheaf on a topological space $X$ defined by $\cals(X) = G$
and $\cals(U)=0$ for all $U \ne X$, then all the stalks $\cals_p$
vanish, but $\cals$ is not the zero presheaf.
\end{exam}

A \emph{sheaf} is a presheaf with two additional properties that
link the global and local sections of the presheaf.
In practice,
most of the presheaves one encounters are sheaves.
Unlike the example in the preceding paragraph, a sheaf all of whose
stalks vanish has no nonzero global sections (Example~\ref{exam:vanish}).

\begin{definition}
A \term{sheaf}\index{sheaf} $\calf$ of abelian groups on a
topological space $X$ is a presheaf satisfying two additional
conditions for any open set $U \subset X$ and any open cover $\{
U_i \}$ of $U$: 
\enum 
\item[(i)] (uniqueness axiom) if $s, t \in \calf(U)$
are sections such that $s|_{U_i} =t|_{U_i}$ for all $i$, then $s = t$; 
\item[(ii)] (gluing axiom) if $\{ s_i \in \calf(U_i)\}$ is a
collection of sections such that 
\[
s_i|_{U_i \cap U_j} = s_j|_{U_i \cap U_j} \text{ for all } i, j,
\]
then there is a section $s \in
\calf(U)$ such that $s|_{U_i} = s_i$ for each $i$.
\eenum
\end{definition}

Suppose there is an ordering on the index set $I$ of the open
cover $\{ U_i\}_{i\in I}$, and consider the sequence of maps
\begin{equation}\label{6e:sheaf}
0 \to \calf(U) \overset{r}{\to} \prod_i \calf(U_i)
\overset{\delta}{\to} \prod_{i < j} \calf(U_i \cap U_j),
\end{equation}
where $r$ is the restriction $r(\omega)_i = \omega|_{U_i}$ and
$\delta$ is the \v{C}ech coboundary operator
\[
(\delta \omega)_{ij} := \omega_j |_{U_{ij}} - \omega_i|_{U_{ij}}.
\]
Then the two sheaf axioms (i) and (ii) are equivalent to the
exactness of the sequence \eqref{6e:sheaf} at $\calf(U)$ and at
$\prod_i \calf(U_i)$, respectively; i.e., the map $r$ is injective
and $\ker \delta = \im r$.

\begin{example*}
For any open subset $U$ of a topological space $X$, let $\calf(U)$
be the abelian group of constant real-valued functions on $U$. If
$V \subset U$, let $\rho_V^U\colon \calf(U) \to \calf(V)$ be the
restriction of functions. Then $\calf$ is a presheaf on $X$. 
Suppose $X$ has nonempty disjoint subsets (for example, 
$X = \R^n$ with the standard topology).  
Then the
presheaf $\calf$ satisfies the uniqueness axiom but not the gluing
axiom of a sheaf: if $U_1$ and $U_2$ are disjoint open sets in
$X$, and $s_1\in \calf(U_1)$ and $s_2\in \calf(U_2)$ have
different values, then there is no constant function $s$ on $U_1
\cup U_2$ that restricts to $s_1$ on $U_1$ and to $s_2$ on $U_2$.
\end{example*}

\begin{example*}
Let $\underline{\R}$ be the presheaf on a topological space $X$
that associates to every open set $U \subset X$ the abelian group
$\underline{\R}(U)$ consisting of all \emph{locally} constant real-valued
functions
on $U$. Then $\underline{\R}$ is a sheaf. More generally, if $G$
is an abelian group, then the presheaf $\underline{G}$ of locally
constant functions with values in $G$ is a sheaf, called the
\term{constant sheaf} with values in $G$.
\end{example*}

\begin{example*}
The zero presheaf $0$ on a topological space in
Example~\ref{exam:zero} is a sheaf.
\end{example*}

\begin{example*}
The presheaf $\cala^k$ on a manifold that assigns to each open set
$U$ the abelian group of $\cinf$ $k$-forms on $U$ is a sheaf.
\end{example*}

\begin{example*}
The presheaf $\calz^k$ on a manifold that associates to each open
set $U$ the abelian group of closed $\cinf$ $k$-forms on $U$ is a
sheaf.
\end{example*}

\section{The Sheaf Associated to a Presheaf} \label{ss:associated}

Associated to a presheaf $\calf$ on a topological space $X$ is another
topological space $E_{\calf}$, called the \term{\'etal\'e space} 
of $\calf$.%
\index{etale space@\'etal\'e space}
As a set, the \'etal\'e space $E_{\calf}$ is the disjoint union
$\coprod_{p\in X} \calf_p$ of all the stalks of $\calf$.
There is a natural projection map 
$\pi\colon E_{\calf} \to X$ that maps $\calf_p$ to $p$.
A \term{section} of the \'etal\'e space $\pi\colon E_{\calf} \to X$
over $U \subset X$ is a map $s\colon U \to E_{\calf}$ such that
$\pi \comp s = \id_U$, the identity map on $U$.
For any open set $U \subset X$, element $s\in \calf(U)$, and point $p
\in U$, let $s_p \in \calf_p$ be the germ of $s$ at $p$.
Then the element $s\in \calf(U)$ defines a section
$\tilde{s}$ of the \'etal\'e
space over $U$,
\begin{align*}
\tilde{s}\colon U &\to E_{\calf},\\
p &\mapsto s_p \in \calf_p.
\end{align*}
The collection
\[
\{ \tilde{s}(U) \mid U \text{ open in } X, \ s \in \calf(U) \}
\]
of subsets of $E_{\calf}$ satisfies the conditions to be a basis
for a topology on $E_{\calf}$.
With this topology, the \'etal\'e space $E_{\calf}$ becomes
a topological space.
By construction, the topological space $E_{\calf}$
is locally homeomorphic to $X$.
For any element $s\in \calf(U)$, the function $\tilde{s}\colon U \to E_{\calf}$
is a continuous section of $E_{\calf}$.
A section $t$ of the \'etal\'e space $E_{\calf}$ is continuous if
and only if every point $p \in X$ has a neighborhood $U$ such that $t
= \tilde{s}$ on $U$ for some $s \in \calf(U)$.

Let $\calf^+$ be the presheaf that associates to each
open subset $U \subset X$ the abelian group
\[
\calf^+(U) := \{ \text{continuous sections } t\colon U \to
E_{\calf} \}.
\]
Under pointwise addition the presheaf $\calf^+$ is easily seen to be a
sheaf, called the
\term{sheafification} or the \term{associated sheaf} of the presheaf
$\calf$.%
\index{shefification}\index{associated sheaf} There is an obvious presheaf morphism $\theta\colon \calf \to
\calf^+$ that sends a section $s\in \calf(U)$ to the section
$\tilde{s} \in \calf^+(U)$.


\begin{example*}
For an open set $U$ in a topological space $X$,
let $\calf(U)$ be the group of all  
\emph{constant} real-valued functions on $U$.
At each point $p \in X$, the stalk $\calf_p$ is $\R$.
The \'etal\'e space $E_{\calf}$ is $X \times \R$,
but not with its usual topology.
A set in $E_{\calf}$ is open if and only if it is of the form
$U \times \{ r\}$ for an open set $U \subset X$ and a 
number $r \in \R$.
Thus, the topology on $E_{\calf}= X \times \R$ is the
product of the given topology on $X$ and the discrete topology on
$\R$. 
The presheaf $\calf$ is not a sheaf.
The sheafification $\calf^{+}$ is the constant sheaf $\underline{\R}$.
\end{example*}

\begin{figure}[!ht]
\centering
\begin{pspicture}(-2,-1)(4.2,2.2)
\psaxes[ticks=none,labels=none,linecolor=gray]{->}(0,0)(-2,-1)(4,2) \psplot[plotstyle=curve,linewidth=.05]{-2}{0}{0}
\psplot[plotstyle=curve,linewidth=.05,xunit=2cm,yunit=2cm]{0.00001}{1.7}{2.718281828 1 neg x div exp} \psline[linestyle=dashed](-2,1.4)(4,1.4)
\uput{.2}[270](4,0){$x$} \uput{.2}[0](0,2){$y$} \uput{.2}[315](0,1.4){$1$}
\end{pspicture}
\caption{A $\cinf$ function all of whose derivatives vanish at $0$}\label{1cinf}
\end{figure}
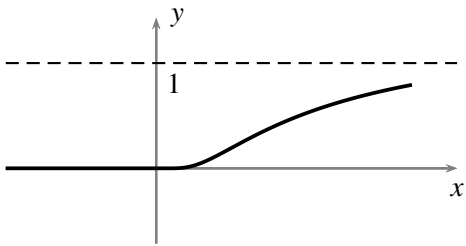

\begin{remark}
The topology of the \'etal\'e space of a sheaf can be quite weird.
For example, if $\calf$ is the sheaf of $\cinf$ real-valued functions
on $\R$, then the \'etal\'e space $E_{\calf}$ is not Hausdroff:
Define $f\colon \R \to \R$ by
\[
f(x)=
\begin{cases}
e^{-1/x}&\quad\text{for } x >0,\\
0 &\quad\text{for } x \le 0.
\end{cases}
\]
Then $f$ is $\cinf$ on $\R$ (see \cite[Example 1.3, p.~5]{tu_m}).
The germ of the function $f$ at $0$ and the germ of the zero
function at $0$ cannot be separated by open sets.
Indeed, if $U$ is any neighborhood of the germ $f_0$ and
$V$ is any neighborjhood of the germ $0_0$, then there exists
an $\epsilon < 0$ such that $f_{\epsilon} = 0_{\epsilon}\in U$ and
$0_{\epsilon} \in V$.
\end{remark}

\begin{exer}
Prove that if $\calf$ is a sheaf,
then $\calf = \calf^+$.
(\textit{Hint}: Show that every continuous section $t\colon U \to
E_{\calf}$
is $\tilde{s}$ for some $s\in \calf(U)$.)
\end{exer}

\prop For every sheaf $\calg$ and every presheaf morphism
$\varphi\colon \calf \to \calg$, there is a unique sheaf morphism
$\varphi^+\colon \calf^+ \to \calg$ such that the diagram
\begin{equation}\label{6e:sheafify}
\xymatrix{
\calf^+ \ar@{-->}[dr]^-{\varphi^+} & \\
\calf \ar[u]^-{\theta} \ar[r]_-{\varphi} & \calg }
\end{equation}
commutes. \eprop

\pf The proof is straightforward and is left as an exercise.
\qed\epf

\section{Sheaf Morphisms}

Recall that all our sheaves are sheaves of abelian groups.
A \term{morphism} $\varphi\colon \calf \to \calg$ of sheaves,
also called a \term{sheaf map}, is by
definition a morphism of presheaves. If $\varphi\colon \calf \to
\calg$ is a morphism of sheaves, then the \term{presheaf kernel}
\[
U \mapsto \ker\big(\varphi_U\colon \calf(U) \to \calg(U)\big)
\]
is a sheaf, called the \term{kernel} of $\varphi$ and written
$\ker \varphi$. The \term{presheaf image}
\[
U \mapsto \im \big(\varphi_U\colon \calf(U) \to \calg(U)\big),
\]
however, is not always a sheaf. The \term{image} of $\varphi$,
denoted $\im \varphi$, is defined to be the sheaf associated to
the presheaf image of $\varphi$.

A sheaf $\calf$ over a space $X$
is a \term{subsheaf} of a sheaf $\calg$ if for every
open set $U$ in $X$ the group $\calf(U)$ is a subgroup
of $\calg(U)$, and the inclusion map $i\colon \calf \to \calg$
is a presheaf morphism.
If $\calf$ is a subsheaf of $\calg$, the \term{quotient sheaf}
is defined to be the sheaf associated to the presheaf
$U \mapsto \calg(U)/\calf(U)$.

A morphism of sheaves $\varphi\colon \calf \to \calg$ is said to
be \term{injective}\index{injective morphism of
sheaves}\index{sheaf morphism!injective} if $\ker \varphi = 0$,
and \term{surjective}\index{surjective morphism of
sheaves}\index{sheaf morphism!surjective} if $\im \varphi =
\calg$.

\prop \label{p:stalk} 
\mbox{}
\enum 
\item[(i)]A morphism of sheaves
$\varphi\colon \calf \to \calg$ is injective if and only if the
stalk map $\varphi_p\colon \calf_p \to \calg_p$ is injective for
every $p$. \item[(ii)] A morphism of sheaves $\varphi\colon \calf
\to \calg$ is surjective if and only if the stalk map
$\varphi_p\colon \calf_p \to \calg_p$ is surjective for every $p$.
\eenum \eprop

\pf Exercise (see \cite[Exercise~1.2 (a),(b), p.\ 66]{hartshorne}).
\qed\epf

In this proposition, neither (i) nor (ii) are true for morphisms of
presheaves.
Let $G$ be an abelian group and $\cals$ the presheaf on a topological
space $X$ such that $\cals(X) =G$ and $\cals(U) = 0$ for all open sets
$U \ne X$ (This is the presheaf of Example~\ref{exam:nonzero}).
Then the stalks $\cals_p$ are all zero.
A counterexample of Proposition~\ref{p:stalk}(i) for presheaves
is $\cals \to 0$; a counterexample to Proposition~\ref{p:stalk}(ii)
for presheaves is $0 \to \cals$.
It is the truth of this proposition for sheaves
that makes sheaves so much more useful than general presheaves.

\begin{exam}\label{exam:vanish}
If the stalk $\calf_p$ of a sheaf $\calf$ vanishes for every $p \in
X$,
then by Proposition~\ref{p:stalk},
the sheaf map $\calf \to 0$ is both injective and surjective,
since the stalk maps $\calf_p \to 0_p= 0$
are injective and surjective for all $p\in X$.
Hence, $\calf$ is isomorphic to the zero
sheaf
and has no nonzero global sections.
\end{exam}

\section{Exact Sequences of Sheaves} \label{s:exact}

A sequence of sheaves and (pre)sheaf morphisms
\[
\cdots \longrightarrow \calf^1 \overset{d_1}{\longrightarrow}
\calf^2 \overset{d_2}{\longrightarrow}  \calf^3
\overset{d_3}{\longrightarrow} \cdots
\]
on a topological space $X$ is said to be \term{exact}\index{exact
sequence!of sheaves} at $\calf^k$ if $\im d_{k-1} = \ker d_k$;
the sequence is said to be \term{exact} if
it is exact at every $\calf^k$. By
Proposition~\ref{p:stalk}, the exactness of a sequence of sheaves
on $X$ is equivalent to the exactness of the sequence of stalk
maps at every point $p \in X$. An exact sequence of sheaves of the
form
\begin{equation} \label{e:short1}
0 \to \cale \to \calf \to \calg \to 0
\end{equation}
is said to be a \term{short exact sequence}.\index{short exact sequence!of
sheaves}
The exactness of a sequence of groups is defined in the same way.

It is not too difficult to show that the exactness of the sheaf
sequence \eqref{e:short1} over a topological space $X$
implies the exactness of the
sequence of sections
\begin{equation} \label{e:left}
0 \to \cale(U) \to \calf(U) \to \calg(U)
\end{equation}
for every open set $U \subset X$,
but the last map $\calf(U) \to \calg(U)$ need not
be surjective.
In fact, as we see in \cite[Theorem~2.8]{elzein--tu},
the cohomology $H^1(U, \cale)$ is a measure of the
nonsurjectivity of the map of global sections $\calf(U) \to
\calg(U)$.

Fix an open subset $U$ of a topological space $X$.
To every sheaf $\calf$ on $X$, we can associate the
abelian group $\Gamma(U,\calf) := \calf(U)$
of sections over $U$
and to every sheaf map $\varphi\colon \calf \to \calg$,
the group homomorphism $\varphi_U\colon \Gamma(U,\calf) \to
\Gamma(U,\calg)$.
This makes $\Gamma(U,\ )$ a functor from sheaves to abelian
groups, henceforth called the \term{section functor}.

\begin{example*}
Let $\calo$ be the sheaf of holomorphic functions on the complex plane
$\C$
and $\calo^*$ the sheaf of nowhere-vanishing holomorphic
functions on $\C$.
For any open set $U \subset \C$,
if $f\in \calo(U)$, then $\exp 2\pi \I f \in \calo^*(U)$.
The kernel of the sheaf map $\exp 2 \pi \I(\ )$
on any open set $U$ consists of the holomorphic
(hence locally constant) integer-valued
functions on $U$.
Hence, there is an exact sequence of sheaves, called
the \term{exponential sequence},
\[
0 \to \underline{\Z} \to \calo \to \calo^* \to 0.
\]
The surjectivity of $\calo \to \calo^*$ follows
from the fact that simply connected neighborhoods form
a basis at each point and if $U$ is simply connected
and $f \in \calo^*(U)$, then $f(U)$ is a simply
connected set in $\C$ not containing the origin,
and hence $\log f$ is defined on $U$.

If $U$ is the punctured plane $\C - \{0\}$, then
the exponential map 
\[
\exp 2 \pi \I (\ )\colon \calo(U) \to \calo^*(U)
\]
is not surjective, since it is not possible to define
its inverse,  $(1/2\pi \I)\log$, on $\C - \{ 0\}$:
for example, $\log z =\log(re^{\I\theta})$ must be defined as $(\log r) + \I\theta$,
but the angle $\theta$ cannot be defined as a continuous
function around a puncture at the origin. 
\end{example*}

A functor $F$ from the category of sheaves on $X$ to the 
category of abelian groups is said to be \term{exact}
if it maps a short exact sequence of sheaves 
\[
0 \to \cale \to \calf \to \calg \to 0
\]
to a short exact sequence of abelian groups
\[
0 \to F(\cale) \to F(\calf) \to F(\calg) \to 0.
\]
If instead one has only the exactness of
\[
0 \to F(\cale) \to F(\calf) \to F(\calg),
\]
then $F$ is said to be a \term{left-exact functor}.
Thus, the section functor $\Gamma(U, \ )$ is left-exact
but not exact.

\section{Resolutions}

Recall that $\underline{\R}$ is the sheaf of locally constant functions
with values in $\R$ and $\cala^k$ is the sheaf of $\cinf$
$k$-forms on a manifold $M$. 
The exterior derivative $d\colon \cala^k(U) \to \cala^{k+1}(U)$,
defined for every open set $U$ in $M$,
induces a morphism of sheaves $d\colon \cala^k \to \cala^{k+1}$.

\prop 
On any manifold $M$ of dimension $n$, the sequence of
sheaves
\begin{equation} \label{6e:dr}
0 \to \underline{\R} \to \cala^0 \overset{d}{\to} \cala^1
\overset{d}{\to} \cdots \overset{d}{\to} \cala^n \to 0
\end{equation}
is exact. \eprop

\pf Exactness at $\cala^0$ is equivalent to the exactness of the
sequence of stalk maps $\underline{\R}_p \to \cala_p^0
\overset{d}{\to} \cala_p^1$ for all $p \in M$. Fix a point $p \in
M$. Suppose $[f] \in \cala_p^0$ is the germ of a $\cinf$ function
$f\colon U \to \R$, where $U$ is a neighborhood of $p$,
such that $d [f] = [0]$ in $\cala_p^1$. Then there is a
neighborhood $V \subset U$ of $p$ on which $df \equiv 0$. Hence,
$f$ is locally constant on $V$ and $[f] \in \underline{\R}_p$.
Conversely, if $[f] \in \underline{\R}_p$, then $d [f] =0$. This
proves the exactness of the sequence \eqref{6e:dr} at $\cala^0$.

Next, suppose $[\omega] \in \cala_p^k$ is the germ of a $\cinf$
$k$-form on some neighborhood of $p$ such that $d [\omega] = 0 \in
\cala_p^{k+1}$. This means there is a neighborhood $V$ of $p$ on
which $d\omega \equiv 0$. By making $V$ smaller, we may assume
that $V$ is contractible. By the Poincar\'e lemma, $\omega$ is
exact on $V$, say $\omega = d\tau$ for some $\tau \in
\cala^{k-1}(V)$. Hence, $[\omega] = d [\tau] $. \qed\epf

In general, an exact sequence of sheaves
\[
0 \to \cala \to \calf^0 \to \calf^1 \to \calf^2 \to \cdots
\]
on a topological space $X$ is called a
\term{resolution}\index{resolution} of the sheaf $\cala$. 
On a complex manifold $M$ of complex dimension $n$, 
the analogue of the Poincar\'e lemma is
the $\bar{\partial}$-Poincar\'e lemma \cite[p.~25]{griffiths--harris}, 
from which it follows that
for each fixed integer $p \ge 0$, the sheaves $\cala^{p,q}$ of $\cinf$
$(p,q)$-forms on $M$ give rise to a resolution of the sheaf 
$\Omega^p$ of
holomorphic $p$-forms on $M$:
\begin{equation} \label{6e:holomorphic}
0 \to {\Omega^p} \to \cala^{p,0} \overset{\bar{\partial}}{\to}
\cala^{p,1} \overset{\bar{\partial}}{\to} \cdots
\overset{\bar{\partial}}{\to} \cala^{p,n} \to 0.
\end{equation}
The cohomology of the \term{Dolbeault complex}
\[
0 \to \cala^{p,0}(M) \overset{\bar{\partial}}{\to}
\cala^{p,1}(M) \overset{\bar{\partial}}{\to} \cdots
\overset{\bar{\partial}}{\to} \cala^{p,n}(M) \to 0
\]
is by definition the \term{Dolbeault cohomology}
of the complex manifold $M$.
(For $(p,q)$-forms on a complex manifold, see
\cite{griffiths--harris}.)

\addtocontents{toc}{\protect\vspace{\baselineskip}}
\section*{\bf Lecture 2.  Sheaf Cohomology}
\medskip

\section{De Rham Cohomology}

To define sheaf cohomology, one might try to imitate the definition 
of de Rham cohomology.
Recall that the de Rham cohomology of a manifold $M$ may be defined as
follows:
\begin{enumerate}
\item[i)] First take a resolution of the sheaf $\underline{\R}$ of locally
constant functions on $M$ by the sheaves of $\cinf$ forms:
\[
0 \to \underline{\R} \to \cala^0 \overset{d}{\to} \cala^1 \overset{d}{\to} \cala^2 \overset{d}{\to} 
\cdots .
\]
\item[ii)] Apply the global section functor:
\[
0\to \Gamma(M, \underline{\R}) \to \cala^0(M) \overset{d}{\to}
 \cala^0(M) \overset{d}{\to}
\]
\item[iii)] Omitting the initial term $\Gamma(M, \underline{\R})$,
take the cohomology of the resulting complex $\cala^{\bullet}(M)$:
\[
H_{\rm{dR}}^k(M) =
\frac{\ker d_k\colon  \cala^k(M) \to \cala^{k+1}(M)}
{\im d_{k-1}\colon \cala^{k-1}(M) \to \cala^k(M)}.
\]
\end{enumerate}

The trouble with this approach is that there could be many
different resolutions of a given sheaf $\calf$,
and the resulting cohomology groups might not be isomorphic.
Fortunately, every sheaf has a unique, canonical resolution
called the \term{Godement canonical resolution}.
Using the Godement canonical resolution, sheaf cohomology
$H^*(X,\calf)$ will be well defined.

\section{The Godement Canonical Functors}

If $\calf$ is a sheaf on a topological space $X$ and $\cale_{\calf}$
its \'etal\'e space, then
\[
\calf (U) = \calf^+(U) = \{ \text{continuous sections of \  } \cale_{\calf} \to X \}.
\]
Let $\calc^0\calf(U) = \{ \text{all (continuous or discontinuous) sections of } \cale_{\calf} \to
X\}$.
Thus, $\calf$ is a subsheaf of $\calc^0\calf$.
If $\calq^1$ is the quotient sheaf, then there is an exact
sequence of sheaves
\begin{equation} \label{s1}
0 \to \calf \to \calc^0\calf \to \calq^1 \to 0.
\end{equation}
Repeat the construction to $\calq^1$, call $\calc^0\calq^1
= \calc^1\calf$ the \term{first Godement functor of
$\calf$}, and let $\calq^2$ be the quotient of $\calc^0\calq^1$ by $\calq^1$:
\begin{equation} \label{s2}
\xymatrix{
0 \ar[r] &\calq^1 \ar[r] &\calc^0\calq^1 \ar[r] \ar@{=}[d] &\calq^2
\ar[r] & 0.\\
& & \calc^1\calf &&
}
\end{equation}
We can splice together the two short exact sequences \eqref{s1} and
\eqref{s2} to
obtain a four-term exact sequence:
\[
\xymatrix{
0 \ar[r] & \calf \ar[r] & \calc^0\calf \ar@{-->}[rr] \ar@{->>}[dr] &&
 C^1\calf \ar[r]&\calq^2
\ar[r] & 0. \\
& & & \calq^1 \ar@{^{(}->}[ur] & &
}
\]
Repeat the construction on $\calq^2$ and so on:
\[
\xymatrix{
0 \ar[r] &\calq^2 \ar[r] &\calc^0\calq^2 \ar[r] \ar@{=}[d] &\calq^3
\ar[r] & 0.\\
& & \calc^2\calf &&
}
\]
By recursion, we obtain a long exact sequence
\[
0 \to \calf \to \calc^0\calf \to \calc^1\calf \to \calc^2\calf \to
\cdots
\]
of sheaves.
This is the \term{Godement canonical resolution} of $\calf$.

The Godement functors are not easy to picture except
possibly $\calc^0\calf$.
An element of $\calc^0\calf$ is a continuous or discontinuous section of
the \'etal\'e space $\cale_{\calf} \to X$ and hence associates
to each $x \in X$ an element of the stalk $\calf_x$.
Therefore, $\calc^0\calf$ is also the direct product $\prod_{x\in X}
\calf_x$.

\section{Sheaf Cohomology} \label{s:cohomology}

Let $\calf$ be a sheaf of abelian groups on a topological space $X$.
To define the sheaf cohomology $H^*(X, \calf)$,
\begin{enumerate}
\item[(i)] Take the Godement canonical resolution of $\calf$:
\[
0 \to \calf \to \calc^0\calf \overset{\delta}{\to} \calc^1\calf
\overset{\delta}{\to} \calc^2\calf \to \cdots .
\]
\item[(ii)] Apply the global section functor $\Gamma(X,\ )$ to
the Godement resolution:
\[
  0 \to \calf(X) \to \calc^0\calf(X) \overset{\delta}{\to}
\calc^1\calf(X) \overset{\delta}{\to} \calc^2\calf(X)
\overset{\delta}{\to} \ldots.
\]
(Recall that $\Gamma(X, \calf) = \calf(X)$.)
\item[(iii)] Take the cohomology of the resulting complex
$\calc^{\bullet}\calf(X)$ of global sections (omitting the initial
term $\calf(X)$):
\[
H^k(X,\calf) = h^k\big( \calc^{\bullet}\calf(X) \big)=
\frac{ \ker \delta\colon \calc^k\calf(X) \to \calc^{k+1}\calf(X)}
{\im \delta\colon \calc^{k-1}\calf(X) \to \calc^{k}\calf(X)}.
\]
\end{enumerate}

\prop
Let $\calf$ be a sheaf of abelian groups on a topological space $X$.
Then the zeroth cohomology $H^0(X, \calf)$ is the
group of global sections:
\[
H^0(X, \calf) = \calf(X).
\]
\eprop

\pf
Because the global section functor is left exact, the exactness 
of the sequence
\[
0 \to \calf \to \calc^0\calf \to \calc^1\calf
\]
implies the exactness of
\[
0 \to \calf(X) \to \calc^0\calf(X) \overset{\delta}{\to}
\calc^1\calf(X).
\]
Hence,
\[
H^0(X,\calf) = \ker \delta = \calf(X).  \tag*{\qedsymbol}
\]
\epf

\section{Flasque Sheaves}

A sheaf $\calf$ on a topological space $X$ is said to be
\term{acyclic} if $H^k(X,\calf) = 0$ for all $k > 0$.
Acyclic sheaves are sheaves with the simplest possible
cohomology groups. They play a special role in sheaf theory.
We will now introduce two types of acyclic sheaves---flasque
sheaves and fine sheaves.

\begin{definition}
A sheaf $\calf$ on a topological space $X$ is said to be 
\term{flasque} if for every open set $U \subset X$, the
restriction map $\rho_U^X\colon \calf(X) \to \calf(U)$ is surjective.
\end{definition}

\begin{example}
A continuous function such as $\sec\colon (-\pi/2, \pi/2) \to \R$
blows up at the endpoints $-\pi/2$ and $\pi/2$ and so cannot be 
extended to a continuous function on $\R$.
This example shows that the sheaf $\calc$ of continuous functions on $\R$
is not flasque, since $\calc(\R) \to \calc\big((-\pi/2, pi/2)\big)$ is not
surjective.
For the same reason, none of the sheaves $\cala^k$ of smooth $k$-forms
on a manifold is flasque.
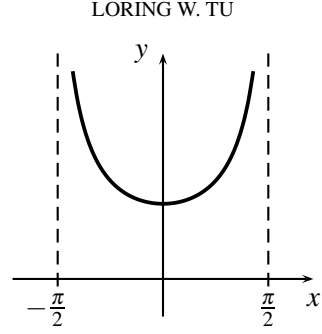
\begin{figure}
\begin{pspicture}(-2,-1)(2,3)
\psline{->}(-2,0)(2,0)
\psline{->}(0,-.5)(0,3)
\psline[linestyle=dashed](-1.4,-.1)(-1.4,3)
\psline[linestyle=dashed](1.4,-.1)(1.4,3)
\psplot[plotstyle=curve,linewidth=.05]%
  {-1.2}{1.2}{1 x 57.32484 mul cos div}
\uput{.2}[270](2,0){$x$}
\uput{.2}[180](0,3){$y$}
\uput{.2}[270](1.4,0){$\frac{\pi}{2}$}
\uput{.2}[270](-1.4,0){$-\frac{\pi}{2}\phantom{-}$}
\end{pspicture}
\caption{The graph of  $y= \sec x$}
\end{figure}
\end{example}

\begin{example}
Recall that for any sheaf $\calf$ on a topological space $X$, we
defined
$\calc^0\calf$ to be the sheaf of all sections of the \'etal\'e space
of $\calf$.
If $U$ is an open set in $X$, then
\[
\calc^0\calf(U) = \prod_{p \in U} \calf_p \quad \text{and} \quad
\calc^0\calf(X) = \prod_{p \in X} \calf_p .
\]
Clearly, the restriction map $\rho\colon \calc^0\calf(X) \to
\calc^0\calf(U)$ is surjective.  
Thus, for any sheaf $\calf$, the sheaf $\calc^0\calf$ is flasque.
\end{example}

Since the $k$th Godement sheaf $\calc^k\calf$ is $\calc^0 \calq^{k-1}$
for some sheaf $\calq^{k-1}$, the Godement sheaf $\calc^k\calf$ is
flasque
for any $k \ge 0$.

We summarize several important properties of flasque sheaves and
the Godement sheaves $\calc^k\calf$ in the proposition below.

\prop
  Let $\calf$ be a sheaf on a topological space $X$.
\begin{enumerate}
\item[(i)] For any $k \ge 0$, the Godement sheaf $\calc^k \calf$ is
  flasque.
\item[(ii)] The Godement functor $\calc^k(\ )$ is an exact functor
  from
sheaves to sheaves \cite[Prop.~2.2.1]{elzein--tu}.
\item[(iii)] The Godement section functor $\Gamma\big(X, \calc^k(\
  )\big)$
is an exact functor from sheaves to abelian groups \cite[Cor.~2.2.7]{elzein--tu}.
\item[(iv)] A flasque sheaf is acyclic \cite[Cor.~2.2.5]{elzein--tu}.
\end{enumerate}
\eprop

\section{Fine Sheaves}

A morphism $f\colon \calf \to \calg$ of sheaves over a  topological
space $X$ induces at each point $x \in X$ a stalk map $f_x\colon
\calf_x \to \calg_x$.
Unlike in the case of continuous real-valued functions, the zero set of
stalk maps
\[
Z= \{ x \in X | f_x = 0\}
\]
is an open subset of $X$.
This is because if two stalk maps agree at a point $p$,
then being germs, they agree in a neighborhood of the point.
Since $Z$ is where the stalk map $f_x$ agrees with the zero map $0_x$,
we see that $Z$ is open in $X$.  Thus,
\[
\supp f = \{ x\in X | f_x \ne 0 \}
\]
is a closed subset of $X$.
(Unlike the support of a real-valued function, the
support of a sheaf morphism is automatically closed without having
to take closure.)

\begin{definition}
Let $\calf$ be a sheaf of abelian groups on a topological space $X$
and $\{ U_{\ga} \}$ a locally finite open cover of $X$.
A \term{partition of unity} for $\calf$ subordinate to $\{ U_{\ga} \}$
is a collection $\{ \eta_{\ga}\colon \calf \to \calf\}$ of sheaf
morphisms such that
\begin{enumerate}
\item[(i)] $\supp \eta_{\ga} \subset U_{\ga}$,
\item[(ii)] at each point $x\in X$, the stalk maps $\eta_{\ga,x}$ sum
  up to the identity map on the stalk $\calf_x$: $\sum \eta_{\ga,x} =
  \id_{\calf_x}$.
\end{enumerate}
\end{definition}

The local finiteness condition guarantees that the sum $\sum_{\ga}
\eta_{\ga,x}$ is a finite sum, since $x$ has a neighborhood that meets
only finitely many of the $U_{\ga}$ and $\supp \rho_{\ga} \subset
U_{\ga}$.

\begin{definition}
A sheaf $\calf$ on a topological space $X$ is \term{fine} if for every
locally finite open cover $\{ U_{\ga}\}$ of $X$, the sheaf $\calf$
admits a partition of unity subordinate to $\{U_{\ga}\}$.
\end{definition}

\prop \label{p:forms}
For each integer $k \ge 0$, the sheaf $\cala^k$ of smooth $k$-forms on
a manifold $M$ is a fine sheaf on $M$.
\eprop

\pf
Let $\{ U_{\ga}\}$ be a locally finite open cover of $M$.
A partition of unity on $M$ subordinate to $\{ U_{\ga}\}$ is a collection of functions
$\{ \rho_{\ga}\colon M \to \R \}$ such that $\supp \rho_{\ga} \subset
U_{\ga}$
and $\sum_{\ga} \rho_{\ga} =1$.
By \cite[Th.~13.7, p.~147]{tu_m}, there is a partition of unity
$\{\rho_{\ga}\}$ on $M$ subordinate to $\{ U_{\ga}\}$.
We will show that a $\cinf$ partition of unity on the manifold $M$
gives rise to a partition of unity for the sheaf $\cala^k$.
Let $\calf = \cala^k$.
Define a collection of morphisms $\{ \eta_{\ga}\colon \calf \to
\calf\}$ for $\calf$ by
\begin{align*}
\eta_{\ga,U}\colon \calf(U) &\to \calf(U), \\
   \eta_{\ga,U}(s) &=
\rho_{\ga} s   \ \text{for } s \in \calf(U).
\end{align*}
Then
\[
\supp \eta_{\ga} \subset \supp \rho_{\ga} \subset U_{\ga}
\]
and
\[
\sum_{\ga} \eta_{\ga,x}(v) = \sum_{\ga} \rho_{\ga} v = v.
\]
Hence, $\{ \eta_{\ga} \}$ is a partition of unity for $\calf$.
\qed\epf

A topological space $X$ is \term{paracompact} if every open cover of
$X$
has a locally finite open refinement.
Continuous partitions of unity exist on a paracompact Hausdorff space
\cite[Th.~41.7]{munkres}.
Thus, if $\cals^k$ is the sheaf of continuous $k$-cochains on a
paracompact Hausdorff space, then the same proof as for 
Proposition~\ref{p:forms} with $\calf=\cals^k$ proves that 
$\cals^k$ is a fine sheaf.

Similarly, the sheaves of continuous forms on
a manifold are fine.  However, in the holomorphic or
algebraic category, partitions of unity generally do not
exist, and the sheaves $\Omega^k$ of holomorphic forms
on a complex manifold and the sheaves $\Omega_{\rm alg}^k$
of algebraic forms on a smooth variety are generally not fine.

The most important property of a fine sheaf is that it is acyclic on 
sufficiently nice spaces.

\begin{thm}
Let $\calf$ be a fine sheaf on a topological space $X$ in which every
open subset is paracompact (a manifold has this property).  Then
$H^k(X, \calf) = 0$ for all $k \ge 1$.
\end{thm}

\bigskip
\addtocontents{toc}{\protect\vspace{\baselineskip}}
\section*{\bf Lecture 3.  Hypercohomology and Spectral Sequences}
\medskip

\section{Cohomology of a Double Complex}

A \term{double complex} is a collection of bigraded abelian groups
$K^{p,q}$ with two commuting differentials
\begin{align*}
d\colon& K^{p,q} \to K^{p,q+1} \text{ in the vertical direction,
  and}\\
\delta\colon& K^{p,q} \to K^{p+1,q} \text{ in the horizontal
  direction},
\end{align*}
i.e., $d^2 = 0$, $\delta^2 = 0$, and $d \delta = \delta d$.
We usually denote a double complex by
\[
K^{\bullet,\bullet} = \bigoplus K^{p,q}
\]
or graphically,
\begin{figure}[!ht]
\begin{center}
\unitlength 1mm
\begin{picture}(55,35)(0,-5)
\multiput(10,0)(10,0){4}{\line(0,1){30}}
\put(5,-5){\makebox(0,0){$0$}}
\put(15,-5){\makebox(0,0){$1$}}
\put(25,-5){\makebox(0,0){$2$}}
\put(35,-5){\makebox(0,0){$3$}}
\put(48,-3){\makebox(0,0){$p$}}
\put(-5,4){\makebox(0,0){$0$}}
\put(-5,12){\makebox(0,0){$1$}}
\put(-5,20){\makebox(0,0){$2$}}
\put(-3,28){\makebox(0,0){$q$}}
\thicklines
\put(0,0){\vector(0,1){30}}
\put(0,0){\vector(1,0){50}}
\put(15,12){\vector(1,0){9}}
\put(15,12){\vector(0,1){9}}
\put(13,16){\makebox(0 ,0){$d$}}
\put(18,10){\makebox(0,0){$\delta$}}
\end{picture}
\caption{}
 \label{8hopf}
\end{center}
\end{figure}

Every row and every column of a double complex $(K^{\bullet,\bullet}, d,\delta)$
is a differential complex. 
The cohomology of the columns is denoted by $H_d$ and the cohomology
of the rows by $H_{\delta}$.
The groups $H_d$ and $H_{\delta}$ each again has a bigrading; for example,
\[
H_d^{p,q} = \frac{\ker d\colon K^{p,q} \to K^{p,q+1}}
{\im d\colon K^{p,q-1} \to K^{p,q}}.
\]

Because $d\delta = \delta d$, the map $\delta$ carries $d$-cocycles to
$d$-cocycles, and $d$-coboundaries to $d$-coboundaries, and induces 
a differential, also denoted by $\delta$, on $H_{d}$.
Thus, $H_{\delta} H_d$ is defined.

Associated to a double complex $K^{\bullet,\bullet}$ is a single
complex $K^{\bullet} = \bigoplus K^k$, where $K^k = \bigoplus_{p+q = k} K^{p,q}$,
obtained by summing along the antidiagonals of $K$.
Because $d$ and $\delta$ commute,
\[
(d + \delta)(d+ \delta) = d\delta + \delta d \ne 0,
\]
so that $\delta + d$ is not a differential.
However, if we alternate the sign of the vertical differential, with 
$D := \delta + (-1)^p d$, then $D^2 = 0$ on the single complex $K^{\bullet}$.
The cohomology $H_D(K^{\bullet})$ of the associated single complex
$K^{\bullet}$ is called the
\term{total cohomology} of the double complex $K^{\bullet,\bullet}$.
Note that the cohomology of the columns is unchanged by switching $d$
for $-d$ every other column.

\section{Cohomology Sheaves}

We have already introduced on several occasions the cohomology
of a complex of abelian groups.  In fact,
in any category in which the notions of kernel and cokernel are
defined,
the cohomology of a complex makes sense.
A \term{complex} in such a category is a sequence of objects and morphisms
\[
C^{\bullet}\colon \cdots \to C^0 \overset{d_0}{\to} C^1
\overset{d_1}{\to} C^2 \to \cdots 
\]
such that the composition of two successive morphisms is zero:
\[
d_k \comp d_{k-1} = 0.
\]
The $k$th cohomology object of the complex $C^{\bullet}$ is
\[
h^k ( C^{\bullet}) = \frac{\ker d_k}{\im d_{k-1}}.
\]
In particular, if $\calc^{\bullet}$ is a complex of sheaves of abelian
groups on a topological space $X$, then we can define the $k$th
\term{cohomology sheaf} of $\calc^{\bullet}$ to be
\[
\calh^k := h^k(\calc^{\bullet})= \dfrac{\ker d_k\colon \calc^k \to
  \calc^{k+1}}{\im d_{k-1}\colon \calc^{k-1} \to \calc^k}.
\]
Note that the \term{cohomology sheaf} $\calh^k$ is a sheaf,
not to be confused with \term{sheaf cohomology} $H^k(X,\calf)$,
which is an abelian group.

\begin{example} \textit{The de Rham complex of sheaves}.
If $\cala^{\bullet}$ is the complex
\[
0 \to \cala^0 \overset{d}{\to} \cala^1\overset{d}{\to} \cala^2 \overset{d}{\to} \cdots
\]
of sheaves of $\cinf$ forms on a manifold $M$, then
by the Poincar\'e lemma the sequence
\[
0\to \underline{\R} \to \cala^0 \overset{d}{\to} \cala^1\overset{d}{\to} \cala^2 \overset{d}{\to} \cdots
\]
is exact.  Therefore, the cohomology sheaves of $\cala^{\bullet}$ are
\[
\calh^k = \begin{cases}
\underline{\R}  &\quad\text{ for } k=0,\\
0 &\quad\text{ for } k>0.
\end{cases}
\]
\end{example}

\section{Filtrations}

A \term{filtration} on an abelian group $A$ is a sequence of
nested subgroups
\[
A = F_0 \supset F_1 \supset F_2 \supset \cdots \supset F_n = 0.
\]
An abelian group with a filtration is called a \term{filtered group}.
Given a filtered group $A$, its \term{associated graded group} is
\[
\Gr(A) = \bigoplus_{p=0}^{n-1} \frac{F_p}{F_{p+1}}.
\]

\begin{figure}
\begin{center}
\begin{pspicture}(0,-2.4)(6,3)
\psline[linewidth=2pt]{->}(0,0)(6,0)
\psline[linewidth=2pt]{->}(0,0)(0,3)
\psline[linewidth=1pt](1,0)(1,3)
\psline[linewidth=1pt](2,0)(2,3)
\psline[linewidth=1pt](3,0)(3,3)
\psline[linewidth=1pt](4,0)(4,3)
\psline[linewidth=1pt](5,0)(5,3)
\uput{.2}[-90](6,0){$p$}
\uput{.2}[180](0,3){$q$}
\uput{.2}[-90](.5,0){$0$}
\uput{.2}[-90](1.5,0){$1$}
\uput{.2}[-90](2.5,0){$2$}
\uput{.2}[-90](3.5,0){$3$}
\psset{linewidth=0.2pt,arrowscale=2,tbarsize=7pt}
  \psline{|<-}(0,-1)(6,-1)\rput*(3,-1){$F_0$}       
  \psline{|<-}(1,-1.7)(6,-1.7)\rput*(3.5,-1.7){$F_1$}       
  \psline{|<-}(2,-2.4)(6,-2.4)\rput*(4,-2.4){$F_2$}       
\end{pspicture}
\caption{The filtration by $p$}
\label{8filtration}
\end{center}
\end{figure}
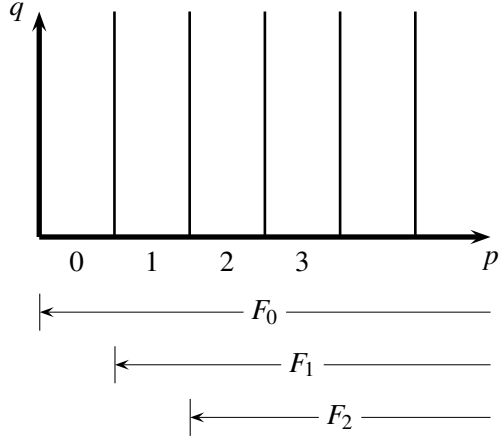

A double complex $\bigoplus K^{p,q}$ has a natural filtration by $p$
as in Figure~\ref{8filtration}.
Each $F_p$ is actually a \term{subcomplex} of $(K^{\bullet}, D)$.
Therefore, the filtration $F_{\bullet}$ by $p$ on $K$ induces a
filtration on $H_D(K)$, which we denote by
$F_{\bullet}\big(H_D(K)\big)$.
It is defined by
\[
F_p\big( H_D(K) \big)\colon \text{ image of }\big( H_D(F_p) \to H_D(K)\big).
\]

\section{Spectral Sequences}

A \term{spectral sequence} is a sequence $\{ (E_r, d_r)
\}_{r=0}^{\infty}$ of differential complexes such that each $E_r$
is the cohomology of its predecessor $(E_{r-1}, d_{r-1})$
for all $r \ge 1$.
We will consider only spectral sequences in which each $E_r$
has a bigrading so that $E_r = \bigoplus_{p,q=0}^{\infty} E_r^{p,q}$.

Fix $(p,q)$ and consider $E_r^{p,q}$ as $r \to \infty$.
If the sequence $E_r^{p,q}$ becomes stationary, i.e., there is an
$r_0$
such that
\[
E_{r_0}^{p,q} = E_{r_0+1}^{p,q} = E_{r_0+2}^{p,q} = \cdots,
\]
then we define $E_{\infty}^{p,q}$ to be the stationary value of
$E_r^{p,q}$.
If for every $(p,q)$, the sequence $E_r^{p,q}$ eventually
becomes stationary as $r$ goes to infinity, then $E_{\infty}$
is defined.
It is important to note that there may not exist an $r_0$
such that $E_{\infty}^{p,q} = E_{r_0}^{p,q}$ for all $p, q$.

One of the main tools for computing the total cohomology
of a double complex is Leray's theorem on the spectral
sequence of a double complex \cite[Th.~14.14]{bott--tu}.

\begin{thm}[Leray]
Given a double complex $(K^{\bullet,\bullet}, d, \delta)$ filtered by
$p$, 
there is a spectral sequence
\[
\{ (E_r, d_r\colon E_r^{p,q} \to E_r^{p+r,q-r+1}) \}
\]
with $(E_0, d_0) = (K^{\bullet,\bullet}, d)$,
$(E_1, d_1) = (H_d, \delta)$, and $E_2 = H_{\delta}H_d$ such that
$E_{\infty}$ is defined and isomorphic to the associated graded group of $H_D(K^{\bullet})$
with the induced filtration by $p$. (In this case, it is customary to
say
that the spectral sequence \term{converges} to $H_D(K^{\bullet}$.)
\end{thm}

Instead of filtering the double complex $K^{\bullet, \bullet}$ by
  $p$,
one can also filter it by $q$.
In this case, by interchanging $p$ and $q$ in Leray's theorem,
one obtains a second spectral sequence with $(E_0, d_0)=
(K^{\bullet,\bullet}, \delta)$, $(E_1, d_1) = (H_{\delta}, d)$,
and $E_2 = H_d H_{\delta}$ that converges to the associated graded
group $H_D(K^{\bullet})$ with the induced filtration by $q$.

The existence of these two spectral sequences, both converging to
$H_D(K^{\bullet})$, is a particularly poweful tool
for proving isomorphism theorems.

\section{Hypercohomology of a Complex of Sheaves}

Let $\calf$ be a sheaf of abelian groups on a topological space $X$.
Recall from Section~\ref{s:cohomology} that the \term{sheaf cohomology} $H^*(X, \calf)$ is defined in
three steps:
\begin{enumerate}
\item[(1)] Take the Godement canonical resolution of $\calf$:
\[
\calc^{\bullet} \calf\colon  0 \to \calf \to \calc^0\calf \to \calc^1
\calf \to \cdots .
\]
\item[(2)] Apply the global section functor $\Gamma(X, \ )$ to
    $\calc^{\bullet}\calf$:
\[
0 \to \calc^0\calf(X) \to \calc^1 \calf(X) \to \cdots 
\]
(Note that the initial term $\calf(X)$ has been dropped.)
\item[(3)] Take the cohomology of $\calc^{\bullet}\calf(X)$:
\[
H^k(X, \calf) = h^k\big( \calc^{\bullet} \calf(X)\big).
\]
\end{enumerate}

Let 
\[
\call^{\bullet} \colon  \call^0 \to \call^1 \to \call^2 \to
\cdots
\]
be a complex of sheaves of abelian groups over a topological space
$X$.
The \term{hypercohomology} of the complex $\call^{\bullet}$
generalizes
the sheaf cohomology of a single sheaf.
The hypercohomology $\H^k(X, \call^{\bullet})$ is also defined in
three steps:
\begin{enumerate}
\item[(1)] Take the Godement canonical resolution of every sheaf
  $\call^q$:
\begin{center}
\begin{tikzpicture}[
	b/.style={	anchor=north,
			align=left}
				]
	\matrix (m) [
	matrix of math nodes,
	nodes in empty cells,
	nodes={
	column 4/.style={anchor=east},
		minimum width=1.5ex,
		outer sep=0pt,
		},
	column sep=2ex,
	row sep=1ex,
	]
	{
	&\hspace{-1.25em}q&&&&&\\
	&&&&&&\\
\call^2 && \calc^0\call^2 && \calc^1\call^2 &&\phantom{\calc^1}\\
	&&&&&&\\
\call^1 && \calc^0\call^1 && \calc^1\call^1 &&\phantom{\calc^1}\\
	&&&&&&\\
\call^0 && \calc^0\call^0 && \calc^1\call^0 &&\phantom{\calc^1}\\
	&&&&&&&&\\
\phantom{\call^0}	&\phantom{q}	&0&&1&&\\
		};
\path[-stealth,>=latex]
			(m-3-1.east) edge (m-3-3.west)
			(m-3-3.east) edge (m-3-5.west)
			(m-3-5.east) edge (m-3-7.west)
			(m-5-1.east) edge (m-5-3.west)
			(m-7-1.east) edge (m-7-3.west)
			(m-5-3.east) edge (m-5-5.west)
			(m-7-3.east) edge (m-7-5.west)
			(m-5-5.east) edge (m-5-7.west)
			(m-7-5.east) edge (m-7-7)
			(m-3-1.north) edge (m-1-1.south)
			(m-3-3.north) edge (m-1-3.south)
			(m-3-5.north) edge (m-1-5.south)
			(m-5-1.north) edge (m-3-1.south)
			(m-5-3.north) edge (m-3-3.south)
			(m-5-5.north) edge (m-3-5.south)
			(m-7-1.north) edge (m-5-1.south)
			(m-7-3.north) edge (m-5-3.south)
			(m-7-5.north) edge (m-5-5.south);
\draw [<->,very thick,>=latex] (m-1-2.north) |- (m-8-8.west);
\node[b] at (m-8-8.west){\hspace{-1em}$p$};

			\end{tikzpicture}
\end{center}



\noindent
The vertical maps $\calc^p \call^q \to \calc^p\call^{q+1}$
are induced from $\call^q \to \call^{q+1}$ by the Godement functors
$\calc^p(\ )$.
This gives rise to a double complex $\bigoplus_{p,q=0}^{\infty} \calc^p \call^q$ of
sheaves.
\item[(2)]  Apply the global section functor to $\bigoplus \calc^p\call^q$
to obtain a double complex $\bigoplus \Gamma(X, \calc^p\call^q)$
of abelian groups.
\item[(3)]  Take the total cohomology of the double complex $\bigoplus \Gamma(X,
\calc^p\call^q)$:
\[
\H^k(X, \call^{\bullet}) = H_D \{ \bigoplus  \Gamma(X, \calc^p\call^q)\}.
\]
\end{enumerate}

\begin{example}
If the complex $\call^{\bullet}$ consists of a single sheaf
\[
0 \to \calf \to 0 \to 0 \to \cdots ,
\]
then the double complex $\bigoplus \Gamma(X, \calc^p\call^q)$ has only
one possibly nonzero row, which occurs in degree $q=0$:

\begin{center}
\begin{tikzpicture}[
	b/.style={	anchor=north,
			align=left}
				]
	\matrix (m) [
	matrix of math nodes,
	nodes in empty cells,
	nodes={
	column 1/.style={anchor=east},
	column 3/.style={anchor=center},
		minimum width=1.5ex,
		outer sep=0pt,
		},
	column sep=.5ex,
	row sep=1ex,
	]
	{
	&\hspace{-1.25em}q&&&&&\\
	&&&&&&&&\\
1 && 0&& 0&&0&&\vphantom{0}\\
	&&&&&&&&\\
0 && \calc^0\calf(X) && \calc^1\calf(X)&&\calc^2 \calf(X) &&\vphantom{\calc^3\calf(X)}\\
	&&&&&&&&&&&&\\
\phantom{\call^0}	&\phantom{q}	&0&&1&&2&&\\
		};
\path[-stealth,>=latex]
			(m-3-3.east) edge (m-3-5.west)
			(m-3-5.east) edge (m-3-7.west)
			(m-3-7.east) edge (m-3-9.west)
			(m-5-3.east) edge (m-5-5.west)
			(m-5-5.east) edge (m-5-7.west)
			(m-5-7.east) edge (m-5-9.west)
			(m-3-3.north) edge (m-1-3.south)
			(m-3-5.north) edge (m-1-5.south)
			(m-3-7.north) edge (m-1-7.south)
			(m-5-3.north) edge (m-3-3.south)
			(m-5-5.north) edge (m-3-5.south)
			(m-5-7.north) edge (m-3-7.south);
\draw [<->,very thick,>=latex] (m-1-2.north) |- (m-6-10.west);
\node[b] at (m-6-10){\hspace{-1.5em}$p$};

			\end{tikzpicture}
\end{center}

In this case the total cohomology of the double complex $\bigoplus
\Gamma(X, \calc^p\call^q)$ is the cohomology of the zeroth row.
Hence,
\[
\H^*(X, \call^{\bullet}) = h^*\big(\calc^{\bullet}\calf(X) \big) = H^*(X, \calf).
\]
It is in this sense that hypercohomology generalizes sheaf cohomology.
\end{example}

\bigskip
\addtocontents{toc}{\protect\vspace{\baselineskip}}
\section*{\bf Lecture 4.  Applications}
\bigskip

In the computation of hypercohomology using spectral sequences,
there are two facts about exact functors that we will repeatedly use:
\begin{enumerate}
\item[(i)] The Godement section functors $\Gamma\big(X, \calc^p(\
  )\big)$
are exact functors from sheaves on $X$ to abelian groups \cite[Cor.~2.2.7]{elzein--tu}.
\item[(ii)] Exact functors commute with cohomology \cite[Prop.~2.2.10]{elzein--tu}:  if $T$
is an exact functor from sheaves to abelian groups and
$\call^{\bullet}$ is a complex of sheaves, then
\[
T\big(\calh (\call^{\bullet})\big) = H\big(T(\call^{\bullet})\big).
\]
\end{enumerate}

\section{The de Rham Theorem}

Let $\cala^k$ be the sheaf of $\cinf$ $k$-forms of dimension $n$.
By the Poincar\'e lemma,
\[
0 \to \underline{\R} \to \cala^0 \overset{d}{\to} \cala^1 \overset{d}{\to}
\cala^2 \to \cdots \to \cala^n \to 0
\]
is an exact sequence.
In particular, the cohomology sheaves of the complex $\cala^{\bullet}$ are
\[
\calh^q =
\begin{cases}
\underline{\R} &\text{for $q=0$,}\\
0 &\text{for $q \ge 1$.}
\end{cases}
\]

We will now compute the hypercohomology $\hh^*(M, \cala^{\bullet})$
using the two spectral sequences of the double complex
\[
E_0 = \bigoplus K^{p,q} = \bigoplus \Gamma(M, \calc^p\cala^q).
\]
Since $\Gamma\big(M, \calc^p(\ )\big)$ is an exact functor
and exact functors commute with cohomology,
\[
E_1 = H_d(E_0) = \bigoplus \Gamma(M,  \calc^p\calh^q)=
\begin{cases}
\bigoplus \Gamma(M, \calc^p \underline{\R}),  &\text{for $q = 0$;}\\
0                       &\text{for $q > 0$.}
\end{cases}
\]

\begin{center}
\begin{pspicture}(-1,-.5)(7,2)
\psline{->}(0,-.1)(7,-.1)
\psline{->}(0,-.1)(0,1.7)
\psline(2,-.1)(2,1.7)
\psline(4,-.1)(4,1.7)
\psline(6,-.1)(6,1.7)
\uput{.2}[180](-0,.75){$E_1 =$}
\rput(1,-0.4){$0$}
\rput(3,-0.4){$1$}
\rput(5,-0.4){$2$}
\rput(1,0.25){$\calc^0\underline{\R}(M)$}
\rput(3,0.25){$\calc^1\underline{\R}(M)$}
\rput(5,0.25){$\calc^2\underline{\R}(M)$}
\rput(1,.75){$0$}
\rput(3,.75){$0$}
\rput(5,.75){$0$}
\rput(1,1.25){$0$}
\rput(3,1.25){$0$}
\rput(5,1.25){$0$}
\uput{.2}[270](7,-.1){$p$}
\uput{.2}[180](0,1.7){$q$}
\rput(7.2,-.1){.}
\end{pspicture}
\end{center}

\noindent
Thus,

\begin{center}
\begin{pspicture}(-1,-.5)(7,2)
\psline{->}(0,-.1)(7,-.1)
\psline{->}(0,-.1)(0,1.7)
\psline(2,-.1)(2,1.7)
\psline(4,-.1)(4,1.7)
\psline(6,-.1)(6,1.7)
\uput{.2}[180](-0,.75){$E_2 = H_{\delta} H_d =$}
\rput(1,-0.4){$0$}
\rput(3,-0.4){$1$}
\rput(5,-0.4){$2$}
\rput(1,0.25){$H^0(M, \underline{\R})$}
\rput(3,0.25){$H^1(M, \underline{\R})$}
\rput(5,0.25){$H^2(M, \underline{\R})$}
\rput(1,.75){$0$}
\rput(3,.75){$0$}
\rput(5,.75){$0$}
\rput(1,1.25){$0$}
\rput(3,1.25){$0$}
\rput(5,1.25){$0$}
\uput{.2}[270](7,-.1){$p$}
\uput{.2}[180](0,1.7){$q$}
\rput(7.2,-.1){.}
\end{pspicture}
\end{center}
Since $d_r$ moves down $r-1$ rows, all differentials $d_r$, $r \ge 2$,
are zero.
It follow that the spectral sequence degenerates at the $E_2$ term
and 
\begin{equation} \label{e:isom1}
\hh^k (M,\cala^{\bullet}) \simeq H^k (M, \underline{\R}).
\end{equation}
There is no extension problem because the associated graded group
of $H^k(M,\underline{\R})$ has only one nonzero term, so that
$H^k(M,\underline{\R}) = \Gr\big( H^k(M, \underline{\R})\big)$.

On the other hand, the second spectral sequence of the double complex
$E_0 = \bigoplus \Gamma(M, \calc^p\cala^q)$ starts with
\[
E_1^{p,q} = H_{\delta}^{p,q} = H^p(M, \cala^q) =
\begin{cases}
\cala^q(M) &\text{for $p= 0$,}\\
0 &\text{for $p > 0$,}
\end{cases}
\]
since $\cala^q$ is a fine sheaf and hence acyclic.
Thus,
\begin{center}
\begin{pspicture}(-1,-.5)(7,2)
\psline{->}(0,-.1)(7,-.1)
\psline{->}(0,-.1)(0,1.7)
\psline(2,-.1)(2,1.7)
\psline(4,-.1)(4,1.7)
\psline(6,-.1)(6,1.7)
\uput{.2}[180](-0,.75){$E_1 = H_{\delta} =$}
\rput(1,-0.4){$0$}
\rput(3,-0.4){$1$}
\rput(5,-0.4){$2$}
\rput(1,0.25){$\cala^0(M)$}
\rput(3,0.25){$0$}
\rput(5,0.25){$0$}
\rput(1,.75){$\cala^1(M)$}
\rput(3,.75){$0$}
\rput(5,.75){$0$}
\rput(1,1.25){$\cala^2(M)$}
\rput(3,1.25){$0$}
\rput(5,1.25){$0$}
\uput{.2}[270](7,-.1){$p$}
\uput{.2}[180](0,1.7){$q$}
\end{pspicture}
\end{center}

\noindent
and

\begin{center}
\begin{pspicture}(-1,-.5)(7,2)
\psline{->}(0,-.1)(7,-.1)
\psline{->}(0,-.1)(0,2)
\psline(2,-.1)(2,2)
\psline(4,-.1)(4,2)
\psline(6,-.1)(6,2)
\uput{.2}[180](-0,.95){$E_2 = H_d H_{\delta} =$}
\rput(1,-0.4){$0$}
\rput(3,-0.4){$1$}
\rput(5,-0.4){$2$}
\rput(1,0.3){$H_{\rm dR}^0(M)$}
\rput(3,0.3){$0$}
\rput(5,0.3){$0$}
\rput(1,.95){$H_{\rm dR}^1(M)$}
\rput(3,.95){$0$}
\rput(5,.95){$0$}
\rput(1,1.6){$H_{\rm dR}^2(M)$}
\rput(3,1.6){$0$}
\rput(5,1.6){$0$}
\uput{.2}[270](7,-.1){$p$}
\uput{.2}[180](0,2){$q$}
\rput(7.2,-.1){.}
\end{pspicture}
\end{center}
This spectral sequence also degenerates at the $E_2$ term and
\begin{equation} \label{e:isom2}
\hh^k(M, \cala^{\bullet}) \simeq H_{\rm dR}^k(M).
\end{equation}
Combining \eqref{e:isom1} and \eqref{e:isom2} gives an 
isomorphism $H_{\rm dR}^*(M) \simeq H^k(M, \underline{\R})$ 
between de Rham cohomology and sheaf cohomology.

\begin{thm} \label{t:dr}
On a manifold $M$, there is a canonical isomorphism
between de Rham cohomology $H_{\rm dR}^*(M)$ and the sheaf cohomology 
$H^*(M,\underline{\R})$ of
the sheaf $\,\underline{\R}$ of locally constant functions with values
in $\R$.
\end{thm}

Instead of the complex $\cala^{\bullet}$ of sheaves of $\cinf$
forms on a manifold $M$, one may also consider on a 
topological space $X$ the sheaf $\cals^q$ of continuous $q$-cochains
with values in $\R$.  There is an exact sequence
\[
0 \to \underline{\R} \to \cals^0 \overset{\delta}{\to}
\cals^1 \overset{\delta}{\to} \cals^2 \overset{\delta}{\to} \cdots
\]
of sheaves.  
By the same computation as for $\hh^*(M, \cala^{\bullet})$ above,
one can show that if all the open subsets of $X$ are paracompact, then
\[
\hh^*(X, \cals^{\bullet}) \simeq H^*(X, \underline{\R})
\]
and
\[
\hh^*(X, \cals^{\bullet}) \simeq H_{\rm sing}^*(X, \R).
\]

\begin{thm} \label{t:singular}
On a topological space $X$ in which all the open subsets are paracompact, there is a canonical isomorphism between
singular cohomology with real coefficients and the 
sheaf cohomology of the sheaf $\underline{\R}$:
\[
H_{\rm sing}^k(X, \R) \simeq H^k(X, \underline{\R}).
\]
\end{thm}

Theorems \ref{t:dr} and \ref{t:singular} together prove the de Rham
theorem:  on a manifold $M$, there is a canonical isomorphism
\[
H_{\rm dR}(M) \simeq H_{\rm sing}^k(M, {\R})
\]
between de Rham cohomology and real singular cohomology.

\section{The de Rham Theorem with Complex Coefficients}

On a complex manifold, instead of differential forms with real coefficients,
we can consider forms
with complex coefficient.
Let $\cala_{\C}^k$ be the sheaf of smooth $k$-forms with complex
coefficients on $M$.
By the Poincar\'e lemma with complex coefficients,
the sequence 
\[
0 \to \underline{\C} \to \cala_{\C}^0 \to \cala_{\C}^1 \to 
\cala_{\C}^2 \to \cdots
\]
is exact.
The same computations as above show that the hypercohomology
of $\cala_{\C}^{\bullet}$ is isomorphic to sheaf cohomology or to de
Rham cohomology:
\begin{align*}
\hh^*(M, \cala_{\C}^{\bullet}) &\simeq H^*(M, \underline{\C}),\\
\hh^*(M, \cala_{\C}^{\bullet}) &\simeq H_{\rm dR}^*(M, \C).
\end{align*}
Combined with Theorem~\ref{t:singular}, these isomorphisms yield the
de Rham theorem with complex
coefficients: 
\[
H_{\rm  dR}^* (M, \C)  \simeq H^*(M, \C)\simeq H_{\rm sing}^*(M, \C).
\]

\section{The Analytic de Rham Theorem}

On a complex manifold $M$, let $\Omega^k$ be the sheaf of
holomorphic $k$-forms and let $\cala^{p,q}$ be the sheaf of $\cinf$
$(p,q)$-forms.
The Poincar\'e lemma has a holomorphic analogue: the sequence
\[
0 \to \underline{\C} \to \Omega^0 \overset{d}{\to} \Omega^1
\overset{d}{\to} \Omega^2 \to \cdots
\]
is exact (for a proof, see \cite[Th.~2.5.1]{elzein--tu}).
It also has a $\bar{\partial}$-analogue:
on a complex manifold $M$ of complex dimension $n$, the
sequence
\[
0\to \Omega^p\to \cala^{p,0} \overset{\bar{\partial}}{\to}
\cala^{p,1} \overset{\bar{\partial}}{\to}
\cala^{p,2} \to \cdots \to \cala^{p,n} \to 0
\]
is exact \cite[p.~25]{griffiths--harris}.

By the holomorphic Poincar\'e lemma, the cohomology sheaves of the
complex $\Omega^{\bullet}$ are
\[
\calh^q(\Omega^{\bullet}) = 
\begin{cases}
\underline{\C} &\text{for $q=0$,}\\
0 &\text{for $q > 0$.}
\end{cases}
\]
With the double complex $E_0 = \bigoplus K^{p,q} = \bigoplus \Gamma(M,
\calc^p\Omega^q)$ filtered by $p$,
\[
E_1 = H_d(E_0) = \bigoplus \Gamma(M, \calc^p \calh^q) 
=\begin{cases}
\bigoplus \Gamma(M, \calc^p \C) &\text{for $q=0$}\\
0 &\text{for $q > 0$,}
\end{cases}
\]
or
\begin{center}
\begin{pspicture}(-1,-.5)(7,2)
\psline{->}(0,-.1)(7,-.1)
\psline{->}(0,-.1)(0,1.7)
\psline(2,-.1)(2,1.7)
\psline(4,-.1)(4,1.7)
\psline(6,-.1)(6,1.7)
\uput{.2}[180](-0,.75){$E_1 =$}
\rput(1,-0.4){$0$}
\rput(3,-0.4){$1$}
\rput(5,-0.4){$2$}
\rput(1,0.25){$\Gamma(M,\calc^0\C)$}
\rput(3,0.25){$\Gamma(M,\calc^1\C)$}
\rput(5,0.25){$\Gamma(M,\calc^2\C)$}
\rput(1,.75){$0$}
\rput(3,.75){$0$}
\rput(5,.75){$0$}
\rput(1,1.25){$0$}
\rput(3,1.25){$0$}
\rput(5,1.25){$0$}
\uput{.2}[270](7,-.1){$p$}
\uput{.2}[180](0,1.7){$q$}
\rput(7.2,-.1){.}
\end{pspicture}
\end{center}

\noindent
Hence,
\begin{center}
\begin{pspicture}(-1,-.5)(7,2)
\psline{->}(0,-.1)(7,-.1)
\psline{->}(0,-.1)(0,1.7)
\psline(2,-.1)(2,1.7)
\psline(4,-.1)(4,1.7)
\psline(6,-.1)(6,1.7)
\uput{.2}[180](-0,.75){$E_2 = H_{\delta} H_d =$}
\rput(1,-0.4){$0$}
\rput(3,-0.4){$1$}
\rput(5,-0.4){$2$}
\rput(1,0.25){$H^0(M, \underline{\C})$}
\rput(3,0.25){$H^1(M, \underline{\C})$}
\rput(5,0.25){$H^2(M, \underline{\C})$}
\rput(1,.75){$0$}
\rput(3,.75){$0$}
\rput(5,.75){$0$}
\rput(1,1.25){$0$}
\rput(3,1.25){$0$}
\rput(5,1.25){$0$}
\uput{.2}[270](7,-.1){$p$}
\uput{.2}[180](0,1.7){$q$}
\rput(7.2,-.1){.}
\end{pspicture}
\end{center}
Therefore, the spectral sequence degenerates at the $E_2$ term
and $E_2 = E_{\infty}$, so that there is a group isomorphism
\[
H^k(M, \underline{\C}) \simeq \hh^k(M, \Omega^{\bullet}).
\]

By an earlier result, there is an isomorphism between singular
cohomology and sheaf cohomology
\[
H_{\rm sing}^k(M,\C) \simeq H^k(M, \underline{\C}).
\]
Hence,
\[
H_{\rm sing}^k(M,\underline{\C}) \simeq \hh^k(M, \Omega^{\bullet}).
\]
This is the analytic de Rham theorem.  
It shows that the singular cohomology with complex coefficients of a
complex
manifold can be computed using only holomorphic forms.

Unlike in the smooth case, however, because the sheaves $\Omega^q$
of holomorphic forms are not acyclic, one cannot conclude that
$H_{\rm sing}^k (M, \C)$ is the cohomology of the holomorphic
de Rham complex $\Omega^{\bullet}(M)$.

\section{Acyclic Resolutions}

A resolution
\[
0 \to \calf \to \call^0 \to \call^1 \to \call^2 \to \cdots
\]
of a sheaf $\calf$ is \term{acyclic} if all the sheaves $\call^q$ are acyclic.
Suppose $0 \to \calf \to \call^{\bullet}$ is an acyclic resolution.
To compute the hypercohomology $\hh^*(X,\call^{\bullet})$, we start
with the double complex
\[
E_0 = \bigoplus E_0^{p,q} = \bigoplus \Gamma(X, \calc^p\call^q),
\]
filtered by $p$.  Then
\begin{align*}
E_1 = H_d(E_0) &= \bigoplus H_d\big( \Gamma(X, \calc^p\call^{\bullet})\big)\\
&=\bigoplus \Gamma (X, \calc^p \calh^q), 
\end{align*}
where $\calh^q$ is the $q$th cohomology sheaf of $\call^{\bullet}$
and the last equality follows from the fact that cohomology 
commutes with the exact functor $\Gamma\big(X,\calc^p(\ )\big)$.

Because $0 \to \calf \to \call^{\bullet}$ is exact,
the cohomology sheaves of $\call^{\bullet}$ are
\[
\calh^q =\begin{cases}
\calf &\text{for $q=0$,}\\
0  &\text{for $q > 0$.}
\end{cases}
\]
Therefore, $E_1$ has only one nonzero row and it is 
row 0:
\begin{center}
\begin{pspicture}(-1,-.5)(7,2)
\psline{->}(0,-.1)(7,-.1)
\psline{->}(0,-.1)(0,1.7)
\psline(2,-.1)(2,1.7)
\psline(4,-.1)(4,1.7)
\psline(6,-.1)(6,1.7)
\uput{.2}[180](-0,.75){$E_1 = H_d =$}
\rput(1,-0.4){$0$}
\rput(3,-0.4){$1$}
\rput(5,-0.4){$2$}
\rput(1,0.25){$\Gamma(X,\calc^0\calf)$}
\rput(3,0.25){$\Gamma(X,\calc^1\calf)$}
\rput(5,0.25){$\Gamma(X,\calc^2\calf)$}
\rput(1,.75){$0$}
\rput(3,.75){$0$}
\rput(5,.75){$0$}
\rput(1,1.25){$0$}
\rput(3,1.25){$0$}
\rput(5,1.25){$0$}
\uput{.2}[270](7,-.1){$p$}
\uput{.2}[180](0,1.7){$q$}
\rput(7.2,-.1){.}
\end{pspicture}
\end{center}

\noindent
It follows that
\begin{center}
\begin{pspicture}(-1,-.5)(7,2)
\psline{->}(0,-.1)(7,-.1)
\psline{->}(0,-.1)(0,1.7)
\psline(2,-.1)(2,1.7)
\psline(4,-.1)(4,1.7)
\psline(6,-.1)(6,1.7)
\uput{.2}[180](-0,.75){$E_2 = H_{\delta} H_d =$}
\rput(1,-0.4){$0$}
\rput(3,-0.4){$1$}
\rput(5,-0.4){$2$}
\rput(1,0.25){$H^0(X, \calf)$}
\rput(3,0.25){$H^1(X, \calf)$}
\rput(5,0.25){$H^2(X, \calf)$}
\rput(1,.75){$0$}
\rput(3,.75){$0$}
\rput(5,.75){$0$}
\rput(1,1.25){$0$}
\rput(3,1.25){$0$}
\rput(5,1.25){$0$}
\uput{.2}[270](7,-.1){$p$}
\uput{.2}[180](0,1.7){$q$}
\rput(7.2,-.1){.}
\end{pspicture}
\end{center}
Because each antidiagonal has only one nonzero box,
there is no extension problem and
\begin{equation} \label{e:sheaf}
\hh^k(X, \call^{\bullet}) \simeq H^k(X, \calf).
\end{equation}

On the other hand, we can also filter
\[
E_0 = \bigoplus E_0^{p,q}= \bigoplus \Gamma(X, \calc^p\call^q),
\]
by $q$.
Note that  the cohomology of the $q$th row is precisely the sheaf
cohomology
$H^*(X, \call^q)$:
\[
E_1 = H_{\delta}(E_0) = \bigoplus H^p(X, \call^q).
\]

Since $\call^q$ is acyclic, the only nonzero column is the
zeroth column
\begin{center}
\begin{pspicture}(-1,-.5)(7,2)
\psline{->}(0,-.1)(7,-.1)
\psline{->}(0,-.1)(0,1.7)
\psline(2,-.1)(2,1.7)
\psline(4,-.1)(4,1.7)
\psline(6,-.1)(6,1.7)
\uput{.2}[180](-0,.75){$E_1 = H_{\delta}(E_0) =$}
\rput(1,-0.4){$0$}
\rput(3,-0.4){$1$}
\rput(5,-0.4){$2$}
\rput(1,0.25){$H^0(X,\call^0)$}
\rput(3,0.25){$0$}
\rput(5,0.25){$0$}
\rput(1,.75){$H^0(X, \call^1)$}
\rput(3,.75){$0$}
\rput(5,.75){$0$}
\rput(1,1.25){$H^0(X, \call^2)$}
\rput(3,1.25){$0$}
\rput(5,1.25){$0$}
\uput{.2}[270](7,-.1){$p$}
\uput{.2}[180](0,1.7){$q$}
\rput(7.2,-.1){.}
\end{pspicture}
\end{center}
Thus,
\[
E_2 = H_d(E_1) = h^*\big(\call^{\bullet}(X)\big).
\]
This spectral sequence also degenerates at the $E_2$ page and
\begin{equation} \label{e:sections}
\hh^k(X, \call^{\bullet}) \simeq E_{\infty}^{0,k} \simeq h^k(\call^{\bullet}(X) \big).
\end{equation}

Combining \eqref{e:sheaf} and \eqref{e:sections}, we have proven the 
following theorem.

\begin{thm} \label{t:acyclic}
Sheaf cohomology can be computed as the cohomology of the complex of
global sections of an acyclic resolution:  let 
\[
0 \to \calf \to \call^0 \to \call^1 \to \call^2 \to \cdots
\]
be an acyclic resolution of a sheaf $\calf$ over a topological space
$X$.
Then
\[
H^k(X, \calf) \simeq h^k \big(\call^{\bullet}(X) \big).
\]
\end{thm}

\begin{example}
Since the Godement resolution $\calc^{\bullet}\calf$ of any sheaf
$\calf$
is flasque and hence acyclic, 
applying Theorem~\ref{t:acyclic} to the Godement
resolution returns
\[
H^k(X, \calf) = h^k \big( \calc^{\bullet}\calf (X) \big),
\]
the definition of sheaf cohomology.
\end{example}

\begin{example}
The complex $\cala^{\bullet}$ of sheaves of $\cinf$ forms on a
manifold
$M$ is a resolution of $\underline{\R}$.
Since the sheaves $\cala^q$ are fine sheaves and hence acyclic,
by Theorem~\ref{t:acyclic} there is an isomorphism between sheaf
cohomology and de Rham cohomology:
\begin{equation} \label{e:dr}
H^k(M,\underline{\R}) \simeq h^k\big( \cala^{\bullet}(M) \big)
= H_{\rm dR}^k(M).
\end{equation}
\end{example}

\begin{example}
The complex $\cals^{\bullet}$ of sheaves of continuous cochains 
on a topological space $X$ all of whose open subsets are paracompact
 is an acyclic resolution of
$\underline{\R}$.
By Theorem~\ref{t:acyclic}
\begin{equation} \label{e:singular}
H^k(X, \underline{\R}) \simeq h^k\big( \cals^{\bullet}(X) \big)
= H_{\rm sing}^k(X, \R).
\end{equation}
Combining \eqref{e:dr} and \eqref{e:singular} gives the de Rham
theorem again.
\end{example}

\begin{example}
In the holomorphic category there do not generally exist
 partitions of unity and the
sheaves
$\Omega^k$ of holomorphic $k$-forms are in general not acyclic.
This explains why the analytic de Rham theorem does not have
the same form as the usual de Rham theorem.
\end{example}

\begin{example}
Let $\cala^{p,q}$ be the sheaf of $\cinf$ $(p,q)$-forms on a complex
manifold $M$ of complex dimension $n$.
By the $\bar{\partial}$-Poincar\'e lemma,
\[
0 \to \Omega^p \to \cala^{p,0} \overset{\bar{\partial}}{\to} \
\cala^{p,1}\overset{\bar{\partial}}{\to} \cdots \to
\cala^{p,n}
\to 0
\]
is a resolution of $\Omega^p$.
Since the $\cala^{p,q}$ are fine sheaves and hence acyclic,
by Theorem \ref{t:acyclic},
\[
H^q(M, \Omega^p) \simeq h^q\left( \cala^{p,\bullet}(M) \right)
 \simeq H^{p,q}(M),
\]
which is Dolbeault's theorem.
\end{example}

\end{document}